\documentclass[a4paper,11pt]{article}

\usepackage{amsmath,amsthm,amssymb}
\usepackage{verbatim}

\newtheorem{theorem}{Theorem}[section]
\newtheorem{example}[theorem]{Example}

\newtheorem{remark}[theorem]{Remark}

\renewcommand{\le}{\leqslant}
\renewcommand{\ge}{\geqslant}
\newcommand{\bs}{\boldsymbol}
\newcommand{\lang}{\left\langle}
\newcommand{\rang}{\right\rangle}
\newcommand{\parfrac}[2]{\frac{\partial{#1}}{\partial{#2}}} 
\newcommand{\R}{\mathbb{R}}
\newcommand{\C}{\mathcal{C}}
\newcommand{\la}{\lambda}
\newcommand{\vh}{\bs{\vec{h}}}
\newcommand{\vhe}{\vh_\epsilon}
\newcommand{\ve}{\bs{v}_\epsilon}
\newcommand{\f}{\bs{f}}
\renewcommand{\H}{\mathcal{H}}
\newcommand{\He}{\H^\epsilon}
\newcommand{\Oo}{\mathcal{O}_0}
\newcommand{\Uo}{\mathcal{U}_0}
\newcommand{\ltrans}{ \la^\pitchfork }

\numberwithin{equation}{section}

\begin{document}

\author{Ulysse Serres
\footnote{Institut \'Elie Cartan de Nancy UMR 7502, Nancy-Universit\'e, CNRS, INRIA, BP 239, F-54506 Vand{\oe}uvre-l\`es-Nancy Cedex, France;
email: \texttt{ulysse.serres@iecn.u-nancy.fr}}
}
\date{}

 \title{Microlocal normal forms for regular fully nonlinear two-dimensional control systems}
 \maketitle

\begin{abstract}
In the present paper we deal with fully nonlinear two-dimensional smooth control systems with scalar input
$\dot{q} = \bs{f}(q,u)$, $q \in M$, $u \in U$,
where $M$ and $U$ are differentiable smooth manifolds of respective dimensions two and one.
For such systems, we provide two microlocal normal forms, i.e., local in the state-input space,
using the fundamental necessary condition of optimality for optimal control problems: the Pontryagin Maximum Principle.
One of these normal forms will be constructed around a regular extremal and the other one will be constructed around an abnormal extremal.
These normal forms, which in both cases are parametrized only by one scalar function of three variables, lead to a nice expression for the control curvature of the system.
This expression shows that the control curvature, a priori defined for normal extremals, can be smoothly extended to abnormals.
\\
\medskip

\noindent{\bf Keywords:} {Control system, control curvature, feedback-equivalence, Pontryagin Maximum Principle.}
\\
\medskip

\noindent{\bf MSC2000:}
34K35; 37C10; 53A55; 93C10; 93C15
\end{abstract}

%
%
\section{Introduction}
In the present paper smooth objects are supposed to be of class $\C^\infty$.

State-feedback classification of control systems
has been studied by numerous authors for the last 40 years.
Antecedents of this theory can be traced to the work of Kronecker (\cite[1890]{Kronecker1890}) in the classification of the singular pencils of matrices
(see \cite{Gantmacher1959} for details on the subject).
Eighty years after Kronecker, Brunovsky (\cite{Brunovsky1970}) used this classification to obtain normal forms of linear controllable systems,
which now bare his name.
Then, the feedback classification problem for control-affine
systems with scalar input was heavily studied in
\cite{Agrachev1998, Jakubczyk1990, Jakubczyk1998, JakubczykRespondek1991, KupkaPomet1991, KupkaPomet1995, Respondek1998}
where the authors also gave list of normal forms.
Finally, in \cite{AgrachevZelenko2007}, A. A. Agrachev and I. Zelenko completely solved the problem of the local 
classification generic control-affine systems on a $n$-dimensional manifold with scalar input for any 
$n\ge 4$ and with two inputs for $n=4$ and $n=5$ by giving a complete set of invariants for these equivalence 
problems.

The present paper deals with the feedback classification of fully nonlinear two-dimensional control systems with
scalar input.
More precisely,
we aim to find some microlocal forms
for nonlinear smooth control systems of the type
\begin{equation}\label{contsyst}
\dot{q} = \f(q,u), \quad q \in M, \quad u \in U,
\end{equation}
where $M$ and $U$ are connected smooth manifolds of respective dimension two and one
under the regularity assumption of strong convexity
\begin{equation}\label{reg_ass_strong}
\parfrac{\bs{f}(q,u)}{u} \wedge \parfrac{^2\bs{f}(q,u)}{u^2} \neq 0, \quad
\forall~(q,u) \in M \times U.
\end{equation}

\medskip
In Section \ref{micloc_forms} we present our main results in Theorem \ref{normal_case} and Theorem \ref{abnormal_case}.
Theorem \ref{normal_case} gives the first microlocal normal for system (\ref{contsyst}.
This normal form is given around a normal extremal.
Theorem \ref{abnormal_case} gives the second microlocal normal which is given around an abnormal extremal.
Those two microlocal normal forms enable us to obtain a nice expression of the control curvature of system (\ref{contsyst}) in a neighborhood of the extremal along which the normalization has been made.
Moreover, in the abnormal case, this expression shows that the control curvature which is a priori only defined for normal extremals, can be smoothly extended to abnormals.

%
%
\section{Preliminaries}\label{Preliminaries}

\subsection{Counting the principal invariants}
Systems of the form (\ref{contsyst}) are considered up to state-feedback equivalence, i.e.,
up to transformations of the form
$(q,u) \to (\phi(q),\psi(q,u))$,
where $\phi$ is a diffeomorphism of $M$ which plays the role of a change of coordinates 
and $\psi$ is a reparametrization of the set $U$ of controls in a way depending on the
state variable $q \in M$.
First of all, let us roughly estimate the number of parameters (invariants) in this equivalence problem.
If the coordinates on the manifold are fixed, a (germ of) control system of type (\ref{contsyst}) is parametrized
by two functions of three variables, and the group of state-feedback transformations
is parametrized by two functions of two variables and one function of three variables. 
Therefore, we can a priori normalize only one function among the two functions defining control system (\ref{contsyst}).
Thus, we expect to have only $2-1=1$
function of three variables
and a certain number of feedback-invariant functions of less than three variables, in the normal forms.

\subsection{Pontryagin Maximum Principle with boundary conditions}
In this section we present a version of the Pontryagin Maximum Principle with boundary conditions (PMP in the sequel)
which will be our main tool in order to obtain microlocal normal forms for system (\ref{contsyst}).
Denote by
$\pi: T^*M \to M$ is the projection of the cotangent bundle to $M$
and by $s$ the canonical Liouville one-form on $T^*M$,
$s_\la=\la \circ \pi_*$, $\la \in T^*M$.
A time-optimal control problem with general boundary conditions takes the form
\begin{eqnarray}
&& \dot{q} = \bs{f}(q,u), \quad q\in M,\quad u\in U   \label{ControlSystemPrelim}\\
&& q(0)\in N_0, \quad q(t_1)\in N_1,       \label{BCControlSystemPrelim}\\
&& t_1 \to \min \quad \textrm{(or $\max$)}, \label{CostControlSystemPrelim}
\end{eqnarray}
where $N_0$ and $N_1$ are given immersed submanifolds of the state space $M$.
Let
$h_u(\lambda) = \lang \lambda,\bs{f}(q,u) \rang$, $\lambda\in T_q^*M$,
be the control dependent Hamiltonian function associated to the control system (\ref{ControlSystemPrelim})
and denote by $\vh_u$ the corresponding Hamiltonian vector field on $T^*M$
(defined by the rule $i_{\vh_u}ds = -d h_u$).
Suppose now that we want to solve the time-optimal problem (\ref{ControlSystemPrelim})$-$(\ref{CostControlSystemPrelim}), then the following holds.
\begin{theorem}[PMP]\label{PMP}
Let an admissible control $u^*(t)$ be time-optimal.
Then, there exists a Lipschitzian curve
$\lambda_t\in T^*M\setminus\{0\}$
such that the following conditions hold for almost all $t\in [0,t_1]$:
\begin{eqnarray}
&& \dot{\lambda}_t = \vh_{u^*(t)}(\lambda_t),                                    \label{PMPthm1} \\
&& h_{u^*(t)}(\lambda_t) = \max_{u\in U} h_u(\lambda_t) = \nu, \quad \nu \in \R, \label{PMPthm3} \\
&& \la_0 \perp T_{\pi(\la_0)}N_0, \quad \la_{t_1} \perp T_{\pi(\la_{t_1})}N_{1}. \label{PMPthm4}
\end{eqnarray}
\end{theorem}
\begin{remark}{\rm
Condition (\ref{PMPthm1}) of PMP says that the solutions of the optimal control problem
(\ref{ControlSystemPrelim})$-$(\ref{CostControlSystemPrelim}) on $M$ are
projections of the solutions of the Hamiltonian system $\dot{\lambda} = \vh_{u^*}(\lambda)$ on $T^*M$.
Moreover, notice that there are two distinct possibilities for condition (\ref{PMPthm3}) of PMP.
If $\nu \neq 0$, then the curve $\la_t$ is called a {\it normal} extremal.
In this case, one can normalize $\lambda_t$ so that $\nu = 1$ (resp. $-1$) in the case of a minimum (resp. maximum) time problem.
If $\nu = 0$, then the curve $\lambda_t$ is called an {\it abnormal} extremal.
}
\end{remark}

\subsection{Curvature of two-dimensional smooth control systems}\label{Sec_courbure}

In this section, we briefly recall some basic facts concerning the curvature of smooth control systems in dimension two.
>From now, we suppose that $M$ and $U$ are connected smooth manifolds of respective dimension two and one.
Let us fix some notations.
We denote by $[\bs{X},\bs{Y}]$ the Lie bracket (or commutator)
$\bs{X} \circ \bs{Y} - \bs{Y} \circ \bs{X}$ of vector fields $\bs{X}$, $\bs{Y} \in \Vec{M}$.
It is again a vector field and in local coordinates on $M$ the Lie bracket reads
$[\bs{X},\bs{Y}](q) = \parfrac{\bs{Y}}{q}\bs{X}(q) - \parfrac{\bs{X}}{q}\bs{Y}(q)$.
If $\bs{X}$ is a smooth vector field on a manifold, we denote by $L_{\bs{X}}$ the Lie derivative along $\bs{X}$.

Denote by 
$h = \max_{u \in U}\lang \la,\bs{f}(q,u) \rang$, 
the Hamiltonian function resulting from the PMP 
by $\H^{\nu}$ the level set $h^{-1}({\nu}) \subset T^*M$, 
and by $\vh^{\nu}$ the Hamiltonian field associated with the restriction of $h^{\nu}$ to $\H^{\nu}$.
Under the regularity assumption (\ref{reg_ass_strong}) 
and the additional assumption
\begin{equation}\label{reg_ass_2}
\bs{f}(q,u) \wedge \parfrac{\bs{f}(q,u)}{u} \neq 0,
\quad
\forall~(q,u) \in M \times U,
\end{equation}
the Hamiltonian function $h$ has constant sign $\epsilon = \pm 1$ and the curve $\He_q = \He \cap T^*_qM$ admits, up to sign and translation, a natural parameter providing us with a vector field
$\bs{v}^{\epsilon}_{q}$ on $\He_q$ and by consequence with a vertical vector field $\ve$ on $\He$ (see e.g. \cite{AAAbook} for details).
The vector field $\ve$ is characterized by the fact that it is, up to sign, the unique vector field on
$\He$ that satisfies
\begin{equation}\label{w''=-ew+bw'}
L^2_{\ve}{s|_{\He}} = -\epsilon{s|_{\He}} + b L_{\ve}{s|_{\He}},
\end{equation}
where $b$ is a smooth function on the level $\He$.
The vector fields $\vhe$ and $\bs{v}_\epsilon$ which are, by definition, feedback-invariant satisfy the nontrivial commutator relation
\begin{equation}\label{[h,[v,h]]=kv}
\Big[ \vhe , \Big[ \ve , \vhe \Big] \Big] = \kappa \ve,
\end{equation}
where the coefficient $\kappa$ is defined to be {\it the control curvature} or simply {\it the curvature} of system (\ref{contsyst}).
\begin{remark}{\rm
The control curvature is by definition  a feedback-invariant of the control system and a function on $\He$ (and not on $M$ as the Gaussian one). Moreover, $\kappa$ is the Gaussian curvature (lifted on $\He$) if the control system defines a Riemannian geodesic problem.
}\end{remark}

%
%
\section{Microlocal normal forms}\label{micloc_forms}

In this section we present two microlocal (i.e. local in the cotangent bundle over the manifold)
normal forms for control systems of type (\ref{contsyst}) under the regularity assumption (\ref{reg_ass_strong}).
Since the feedback-invariants of such a system are
functions on a three-dimensional bundle over the manifold $M$,  the microlocalization of the problem is clearly
reasonable. Actually, under the considered genericity assumption we may not expect better normal forms.
These two normal forms will enable us to get a nice expression for the curvature in restriction to the extremal
along which the normalization is done.

\subsection{Normal case}
Let $\pi:T^*M\to M$ denote the canonical projection.
Fix a pair $(q_0,u_0) \in M \times U$ and assume that both relations (\ref{reg_ass_strong}) and (\ref{reg_ass_2}) are satisfied
at $(q_0,u_0)$.
Let $\la_0 \in T^*_{q_0}M \cap \He$ be a covector satisfying $\lang \la_0, \f(q_0,u_0) \rang = 0$.
For $\tau$ small enough, define the curve
\begin{eqnarray*}
\ltrans:\tau \mapsto e^{\tau \big[\ve,\vhe\big]}(\la_0) \in \He,
\end{eqnarray*}
where $e^{\tau \big[\ve,\vhe\big]}$ denotes the flow of the commutator of the fields $\ve$ and $\vhe$ which are defined according to Section \ref{Sec_courbure}.
The image $N_0$ of $\pi \circ \ltrans$ is canonically defined on the manifold $M$ and according to (\ref{reg_ass_2}) is transverse to projections onto $M$ of the integral curves of $\vhe$.
We will use the curve $N_0$ in order to define the horizontal axis (with origin at $q_0$) of our system of microlocal coordinates on $M$.
Then, vertical lines will be defined as the time-optimal paths that connect points in $M$ to $N_0$.
In other words, vertical lines are the projection onto $M$ of the extremals of the following time-optimal control problem:
\begin{eqnarray*}
&& \dot{q}=\bs{f}(q,u),\quad q \in M,\quad u\in U,\\
&& q(0)\in N_0,\quad q(t_1)=q_1 \quad \textrm{fixed}, \\ 
&&t_1 \to \min.
\end{eqnarray*}
Let $\la_{x_1} \in T^*_{\pi\circ\ltrans(x_1)}M \cap \He$ be such that $\langle \la_{x_1} , \frac{d}{dx_1}(\pi\circ\ltrans) \rangle = 0$.
For $x_1$, $x_2$ small, define the map $\phi$ by 
\begin{equation}\label{Phiq1q2}
\phi(x_1,x_2) = \pi \circ e^{x_2 \vh}\left(\la_{x_1}\right).
\end{equation}
It follows from (\ref{reg_ass_2}) that the differential $D_{(0,0)}\phi$ is bijective which implies that $\phi$ defines a system of local coordinates in a neighborhood of $q_0$.
Denote by $\Oo$ the preimage of this neighborhood by $\phi$.
In the local coordinates system $(x_1,x_2)$ defined by $\phi$ control system (\ref{contsyst}) reads:
\begin{eqnarray*}
\dot{x}_1 &=& f_1(x_1,x_2,u) \\
\dot{x}_2 &=& f_2(x_1,x_2,u), \quad (x_1,x_2) \in \Oo,
\end{eqnarray*}
where, according to (\ref{Phiq1q2}), $f_1$, $f_2$ satisfy
\begin{equation}\label{f2qu0=1}
f_1(x_1,x_2,u_0)=0, \quad
\parfrac{f_1}{u}(0,0,u_0)=1, \quad
f_2(x_1,x_2,u_0)=1, \quad
\parfrac{f_2}{u}(0,0,u_0)=0.
\end{equation}
Since $\parfrac{f_1}{u}(0,0,u_0)\neq 0$, the feedback transformation
$(x_1,x_2,u)\mapsto\tilde{u}=f_1(x_1,x_2,u)$ is well-defined in a neighborhood of $(0,0,u_0)$ and it brings the system to
\begin{eqnarray}\label{ContSystCoord}
\begin{array}{rcl}
\dot{x}_1 &=& \tilde{u} \\
\dot{x}_2 &=& \tilde{f}_2\left(x_1,x_2,\tilde{u}\right).
\end{array}
\end{eqnarray}
According to the third equality in (\ref{f2qu0=1}), $\tilde{f}_2$ satisfies
$\tilde{f}_2(0,0,0)=1$,
which shows that the function $\tilde{f_2}$ can be written in the form
\begin{equation}\label{tildef2=1-psi}
\tilde{f_2}(x_1,x_2,u)=1-\psi(x_1,x_2,u)u.
\end{equation}
Let $(p,x)=(p_1,p_2,x_1,x_2)$ be a canonical coordinates on $T^*\R^2$.
Taking into account (\ref{tildef2=1-psi}), the control dependent Hamiltonian function for the control system
(\ref{ContSystCoord}) reads
\begin{equation*}
h_u(p,x)=p_1u+p_2(1-\psi(x,u)u).
\end{equation*}
We now prove that the function $\psi(x,u)$ satisfies $\psi(x,0)=0$.
By construction, solutions of the time-optimal control problem
\begin{eqnarray*}
&& \dot{x}_1 = u \\
&& \dot{x}_2 = 1-\psi(x_1,x_2,u)u, \\
&& x(0) \in \R\times\{0\}, \quad x(t_1) \in \R\times\{t_1\},\\
&& t_1 \to \min,
\end{eqnarray*}
is the set of all segment lines included in $\Oo$.
Applying Theorem \ref{PMP} to the above time-optimal problem implies that, along the extremal corresponding to the optimal control $u=0$, the covector $p(t)$ is solution to
\begin{eqnarray}
&& \begin{array}{rcl}
   \dot{p}_1 &=& \displaystyle{-\parfrac{h_{u=0}}{q_1}=0} \\
   \dot{p}_2 &=& \displaystyle{-\parfrac{h_{u=0}}{q_2}=0},
   \end{array} \nonumber \\ 
&& p(0) \in \{0\}\times\R, \quad p(t_1) \in \{0\}\times\R. \label{BoundCondCoord}
\end{eqnarray}
Taking into account that the covector $p(t)$ never vanishes and because $t_1$ is arbitrary, one infers from (\ref{BoundCondCoord}) that the covector corresponding to the optimal control $u=0$ is
\begin{equation}\label{p(t)=(0,1)}
p(t) = \left(0,p_2(t)\right),\quad p_2(t) \neq 0 \quad \forall~t.
\end{equation}
Equation (\ref{p(t)=(0,1)}) implies in particular that the maximality condition $\parfrac{h_{u}}{u}|_{u=0}=0$ is equivalent to
$\psi(x,0)=0$, forall $x \in \Oo$, 
from which it follows immediately that the function $\psi$ can be written $\psi(x,u)=\varphi(x,u)u$.
We now prove that the function $\varphi(x,u)$ never vanishes.
>From the regularity assumption (\ref{reg_ass_2}), it follows that $\bs{f}$ has to satisfy
\begin{equation}\label{d2F=-alpha}
\parfrac{^2{\bs{f}}}{u^2} = -\epsilon\alpha{\bs{f}}-\beta\parfrac{{\bs{f}}}{u},
\end{equation}
where $\alpha=\alpha(x,u)$ is positive.
Equation (\ref{d2F=-alpha}) implies in particular that
$\det(\parfrac{^2{\phi_*\bs{f}}}{u^2},\parfrac{\phi_*\bs{f}}{u})|_{u=0}=-\epsilon\alpha \det({\phi_*\bs{f}},\parfrac{\phi_*\bs{f}}{u})|_{u=0}$,
or equivalently, that $2\varphi(x,0)=\epsilon\alpha(x,0)$, which proves that the function $\varphi(x,u)$ never vanishes (at least in a small enough neighborhood $\Oo\times\Uo$ of zero). We can thus set $\varphi=e^{2a}$, with $a \in \C^\infty(\Oo\times\Uo)$.
Summing up, we have proved the following theorem.
\begin{theorem}\label{normal_case}
Under the regularity assumptions (\ref{reg_ass_strong}) and (\ref{reg_ass_2}) control system (\ref{contsyst}) can be put into
the microlocal normal form
\begin{eqnarray*}
\dot{q}_1 &=& u \\
\dot{q}_2 &=& 1 - \epsilon e^{2a(q_1,q_2,u)}u^2,
\end{eqnarray*}
where $\epsilon=1$ (resp. $\epsilon=-1$) if the curves of admissible velocities of system (\ref{contsyst}) are convex (resp. concave).
\end{theorem}
The curvature of the control system in the normal form (\ref{normal_case}) is also easily computed according to formula
(\ref{[h,[v,h]]=kv}) which leads to
\begin{equation*}
\kappa(q_1,q_2,u) = -\parfrac{^2a}{q_2^2}(q,0) - \left(\parfrac{a}{q_2}(q,0)\right)^2 + O(u).
\end{equation*}
\begin{example}
{\rm
Consider the control system
\begin{eqnarray*}
\dot{q}_1 &=& u \\
\dot{q}_2 &=& 1-e^{a(q_1,q_2)}u^2,\quad u \in \R.
\end{eqnarray*}
This system is just the particular case of the normal form (\ref{normal_case}) when the function $a$ only 
depends on the base point $q \in M$.
The curvature of this system is
\begin{equation}\label{PolK}
\kappa(q_1,q_2,u) =
-\parfrac{^2a}{q_2^2} - \left(\parfrac{a}{q_2}\right)^2
-3e^{2a}\parfrac{^2a}{q_2^2}u^2
-e^{2a}\parfrac{^2a}{q_1\partial q_2}u^3.
\end{equation}
It turns out that, if we ask the curvature to be constant then, this system is feedback-equivalent to the normal form
\begin{eqnarray*}
\dot{q}_1 &=& u \\
\dot{q}_2 &=& 1 - e^{2q_2\sqrt{-\kappa}+g(q_1)}u^2,\quad u \in \R,\quad \kappa \le 0,
\end{eqnarray*}
which is easily obtained asking for the vanishing of non zero degree coefficients in polynomial (\ref{PolK}).
}
\end{example}

\subsection{Abnormal case}
The construction of the microlocal normal form around a regular extremal can easily be adapted in order to get a micro
local form around an abnormal extremal, that is, around an extremal along which the Hamiltonian function of PMP vanishes
identically.
Set $\mathcal{T}M_{\mathrm{ab}}=\{\bs{f}(q,u) \mid \bs{f}(q,u) \wedge \parfrac{\bs{f}}{u}(q,u)=0\}$.
To insure the existence of an abnormal trajectory,
we assume that
$\mathcal{T}M_{\mathrm{ab}}$ defines a codimension one submanifold of $TM$.
We do not repeat the detailed construction but only cite the following theorem.
\begin{theorem}\label{abnormal_case}
Suppose that the regularity assumption (\ref{reg_ass_strong}) holds in a neighborhood of $(q_0,u_0)$ for which $\bs{f}(q_0,u_0) \in \mathcal{T}M_{\mathrm{ab}}$.
Then, control system (\ref{contsyst}) can be put into the microlocal normal form
\begin{equation}\label{NormalFormAb}
\begin{array}{rcl}
\dot{q}_1 &=& u \\
\dot{q}_2 &=& e^{2a(q_1,q_2,u)}(1-u)^2.
\end{array}
\end{equation}
%
\end{theorem}
The curvature of the control system in the normal form (\ref{NormalFormAb}) is also easily computed
according to formula (\ref{[h,[v,h]]=kv}) which leads to
\begin{equation*}
\kappa(q_1,q_2,u) = -\parfrac{^2a}{q_1^2}(q,1) - \left(\parfrac{a}{q_1}(q,1)\right)^2 + O(u-1),
\end{equation*}
which shows in particular that the value $\kappa(q,1)$ is well defined so that the curvature can be smoothly extended along the abnormal trajectory.


\end{document}